\documentstyle[12pt,amscd]{amsart}
\newtheorem{thm}[equation]{Theorem}
\numberwithin{equation}{section}
\newtheorem{cor}[equation]{Corollary}

\newtheorem{rmk}[equation]{Remark}

\newtheorem{lem}[equation]{Lemma}
\newtheorem{conj}[equation]{Conjecture}
\newtheorem{defin}[equation]{Definition}

\newtheorem{prop}[equation]{Proposition}

\begin{document}
\raggedbottom \voffset=-.7truein \hoffset=0truein \vsize=8truein
\hsize=6truein \textheight=8truein \textwidth=6truein
\baselineskip=18truept
\def\mapright#1{\ \smash{\mathop{\longrightarrow}\limits^{#1}}\ }
\def\mapleft#1{\smash{\mathop{\longleftarrow}\limits^{#1}}}
\def\mapup#1{\Big\uparrow\rlap{$\vcenter {\hbox {$#1$}}$}}
\def\mapdown#1{\Big\downarrow\rlap{$\vcenter {\hbox {$\ssize{#1}$}}$}}
\def\mapne#1{\nearrow\rlap{$\vcenter {\hbox {$#1$}}$}}
\def\mapse#1{\searrow\rlap{$\vcenter {\hbox {$\ssize{#1}$}}$}}
\def\mapr#1{\smash{\mathop{\rightarrow}\limits^{#1}}}
\def\ss{\smallskip}
\def\vp{v_1^{-1}\pi}
\def\at{{\widetilde\alpha}}
\def\sm{\wedge}
\def\la{\langle}
\def\ra{\rangle}
\def\on{\operatorname}
\def\spin{\on{Spin}}
\def\kbar{{\overline k}}
\def\qed{\quad\rule{8pt}{8pt}\bigskip}
\def\ssize{\scriptstyle}
\def\a{\alpha}
\def\bz{{\Bbb Z}}
\def\im{\on{im}}
\def\ct{\widetilde{C}}
\def\ext{\on{Ext}}
\def\sq{\on{Sq}}
\def\eps{\epsilon}
\def\ar#1{\stackrel {#1}{\rightarrow}}
\def\br{{\bold R}}
\def\bc{{\bold C}}
\def\bh{{\bold H}}
\def\si{\sigma}
\def\Ebar{{\overline E}}
\def\Sum{\sum}
\def\tfrac{\textstyle\frac}
\def\tb{\textstyle\binom}
\def\Si{\Sigma}
\def\w{\wedge}
\def\equ{\begin{equation}}
\def\b{\beta}
\def\G{\Gamma}
\def\g{\gamma}
\def\psit{\widetilde{\Psi}}
\def\tht{\widetilde{\Theta}}
\def\psiu{{\underline{\Psi}}}
\def\thu{{\underline{\Theta}}}
\def\aee{A_{\text{ee}}}
\def\aeo{A_{\text{eo}}}
\def\aoo{A_{\text{oo}}}
\def\aoe{A_{\text{oe}}}
\def\fbar{{\overline f}}
\def\endeq{\end{equation}}
\def\sn{S^{2n+1}}
\def\zp{\bold Z_p}
\def\A{{\cal A}}
\def\P{{\cal P}}
\def\cj{{\cal J}}
\def\zt{{\bold Z}_2}
\def\bs{{\bold s}}
\def\bof{{\bold f}}
\def\bq{{\bold Q}}
\def\be{{\bold e}}
\def\Hom{\on{Hom}}
\def\ker{\on{ker}}
\def\coker{\on{coker}}
\def\da{\downarrow}
\def\colim{\operatornamewithlimits{colim}}
\def\zphat{\bz_2^\wedge}
\def\io{\iota}
\def\Om{\Omega}
\def\u{{\cal U}}
\def\e{{\cal E}}
\def\exp{\on{exp}}
\def\wbar{{\overline w}}
\def\rbar{{\overline r}}
\def\xbar{{\overline x}}
\def\ybar{{\overline y}}
\def\zbar{{\overline z}}
\def\ebar{{\overline e}}
\def\et{{\widetilde E}}
\def\ni{\noindent}
\def\coef{\on{coef}}
\def\den{\on{den}}
\def\lcm{\on{l.c.m.}}
\def\vi{v_1^{-1}}
\def\ot{\otimes}
\def\psibar{{\overline\psi}}
\def\mhat{{\hat m}}
\def\exc{\on{exc}}
\def\ms{\medskip}
\def\ehat{{\hat e}}
\def\etao{{\eta_{\text{od}}}}
\def\etae{{\eta_{\text{ev}}}}
\def\dirlim{\operatornamewithlimits{dirlim}}
\def\gt{\widetilde{L}}
\def\lt{\widetilde{\lambda}}
\def\sgd{\on{sgd}}
\def\ord{\on{ord}}
\def\gd{{\on{gd}}}
\def\rk{{{\on{rk}}_2}}
\def\nbar{{\overline{n}}}

\def\N{{\Bbb N}}
\def\Z{{\Bbb Z}}
\def\Q{{\Bbb Q}}
\def\R{{\Bbb R}}
\def\C{{\Bbb C}}
\def\l{\left}
\def\r{\right}
\def\ls{\leq}
\def\gs{\geq}
\def\bg{\bigg}
\def\({\bigg(}
\def\[{\bigg[}
\def\){\bigg)}
\def\]{\bigg]}
\def\colon{{:}\;}
\def\lg{\langle}
\def\rg{\rangle}
\def\t{\text}
\def\f{\frac}
\def\mo{\on{mod}}
\def\al{\alpha}
\def\ve{\varepsilon}
\def\bi{\binom}
\def\eq{\equiv}
\def\cs{\cdots}
\def\Remark{\noindent{\it  Remark}}
\hbox{J. Pure Appl. Algebra. 209(2007), no.\,1, 57--69.}
\bigskip
\title[A number-theoretic approach to homotopy exponents of ${\rm SU}(n)$]
{A number-theoretic approach to homotopy exponents of ${\rm
SU}(n)$}
\author{Donald M. Davis}
\address{Department of Mathematics, Lehigh University\\Bethlehem, PA 18015, USA}
\email{dmd1@@lehigh.edu}
\author{Zhi-Wei Sun}
\address{Department of Mathematics, Nanjing University
\\Nanjing 210093, People's Republic of China}
\email{zwsun@@nju.edu.cn}

\thanks {
\newline\indent
The second author is supported by the National Science Fund for
Distinguished Young Scholars (no. 10425103) in P.~R.~China.}

\maketitle
\begin{abstract} We use methods of combinatorial number
theory to prove that, for each $n\ge2$ and any prime $p$, some
homotopy group $\pi_i({\rm SU}(n))$ contains an element of order
$p^{n-1+\ord_p(\lfloor n/p\rfloor!)}$, where $\ord_p(m)$
denotes the largest integer $\al$ such that $p^{\al}\mid m$.
\end{abstract}

\section{Introduction}\label{intro}
Let $p$ be a prime number. The {\it homotopy $p$-exponent} of a
topological space $X$, denoted by $\exp_p(X)$, is defined to be
the largest $e\in\N=\{0,1,2,\ldots\}$ such that some homotopy
group $\pi_i(X)$ has an element of order $p^e$. This concept has
been studied by various topologists (cf. \cite{IMJ}, \cite{Gray},
\cite{Selick}, \cite{CMN1}, \cite{CMN2}, \cite{large},
\cite{Neis},  \cite{Th1}, and \cite{Th2}). The most celebrated
result about homotopy exponents (proved by Cohen, Moore, and
Neisendorfer in \cite{CMN1}) states that $\exp_p(S^{2n+1})=n$ if
$p\not=2$.

The {\it special unitary group} ${\rm SU}(n)$ (of degree $n$) is
the space of all $n\times n$ unitary matrices (the conjugate
transpose of such a complex matrix equals its inverse) with
determinant one. (See, e.g., \cite[p.\,68]{Hu}.) It plays a
central role in many areas of mathematics and physics. The famous
Bott Periodicity Theorem (\cite{Bott}) describes $\pi_i({\rm
SU}(n))$ with $i<2n$. In this paper, we provide a strong and
elegant lower bound for the homotopy $p$-exponent of ${\rm
SU}(n)$.

As in number theory, the integral part of a real number $c$ is
denoted by $\lfloor c\rfloor$. For a prime $p$ and an integer $m$,
the {\it $p$-adic order} of $m$ is given by
$\ord_p(m)=\sup\{n\in\N:\,p^n\mid m\}$ (whence
$\ord_p(0)=+\infty$).

Here is our main result.
\begin{thm} For any prime $p$ and $n=2,3,\ldots$,
some homotopy group $\pi_i({\rm SU}(n))$ contains an element of
order $p^{n-1+\ord_p(\lfloor n/p\rfloor!)}$; i.e., we have the
inequality
$$\exp_p({\rm SU}(n))\ge n-1+\ord_p\l(\l\lfloor\f np\r\rfloor!\r).$$
 \label{mainthm}\end{thm}
\noindent We discuss in Section \ref{outlinesec} the extent to
which Theorem \ref{mainthm} might be sharp.

Our reduction from homotopy theory to number theory involves
Stirling numbers of the second kind. For $n,k\in\N$ with
$n+k\in\Z^+=\{1,2,3,\ldots\}$, the Stirling number $S(n,k)$ of the
second kind is the number of partitions of a set of cardinality
$n$ into $k$ nonempty subsets; in addition, we define $S(0,0)=1$.
We will use the following definition.

\begin{defin} Let $p$ be a prime. For $k,n\in\Z^+$ with $k\ge n$,
we define $$e_p(n,k)=\min_{m\ge n}\ord_p(m!S(k,m)).$$\label{edef}
\end{defin}
In Sections \ref{outlinesec} and \ref{p=2sec} we prove the
following standard result.

\begin{prop} \label{redn} Let $p$ be a prime, and let $n\in\Z^+$.
 Then, for all $k\ge n$, we have
$\exp_p({\rm SU}(n))\ge e_p(n,k)$ unless $p=2$ and $n\eq 0\ (\mo\
2)$, in which case $\exp_2({\rm SU}(n))\ge e_2(n,k)-1$.\end{prop}

Our innovation is to extend previous work (\cite{S05}) of the
second author in combinatorial number theory to prove the
following result, which, together with Proposition \ref{redn},
immediately implies Theorem \ref{mainthm} when $p$ or $n$ is odd.
In Section \ref{p=2sec}, we explain the extra ingredient required
to deduce Theorem \ref{mainthm} from \ref{redn} and \ref{numththm}
when $p=2$ and $n$ is even.
\begin{thm}\label{numththm} Let $p$ be any prime and $n$ be a positive integer.

{\rm (i)} For any $\al,h,l,m\in\N$, we have
\begin{eqnarray*}&&\ord_p\(m!\sum_{k=0}^l\bi lk(-1)^kS\l(kh(p-1)p^{\al}+n-1,\,m\r)\)
\\&&\ \ \ge\min\l\{l(\al+1),\,n-1+\ord_p\l(\l\lfloor\f mp\r\rfloor!\r)\r\}.
\end{eqnarray*}

{\rm (ii)} If we define $N=n-1+\lfloor n/(p(p-1))\rfloor$, then
$$e_p\l(n,(p-1)p^L+n-1\r)\ge n-1+\ord_p\l(\l\lfloor\f np\r\rfloor!\r)
\ \ \t{for}\ L=N,N+1,\ldots.$$
\end{thm}

In Section \ref{numthsec}, we prove the following broad
generalization of Theorem \ref{numththm}, and in Section
\ref{outlinesec}, we show that it implies Theorem \ref{numththm}.
\begin{thm}\label{combthm} Let $p$ be a prime, $\al,n\in\N$ and
$r\in\Z$. Then for any polynomial $f(x)\in\Z[x]$ we have
$$\ord_p\(\sum_{k\eq r\,(\mo\ p^{\al})}\bi nk(-1)^kf\l(\f{k-r}{p^{\al}}\r)\)
\gs\ord_p\l(\l\lfloor\f n{p^{\al}}\r\rfloor!\r).$$\end{thm}
\noindent Here we adopt the standard convention that $\bi nk$ is
$0$ if $k$ is a negative integer.

In Theorem \ref{strthm}, we give a strengthened version of Theorem
\ref{combthm}, which we conjecture to be optimal in a certain
sense. Our application to topology uses the case $r=0$ of Theorem
\ref{combthm}; the more technical Theorem \ref{strthm} yields no
improvement in this case.

In \cite{large}, the first author used totally different, and much
more complicated,
 methods
to prove that
\begin{equation}\exp_p({\rm SU}(n))\ge
n-1+\left\lfloor\frac{n+2p-3}{p^2}\right\rfloor
+\left\lfloor\frac{n+p^2-p-1}{p^3}\right\rfloor,
\label{largethm}\end{equation} where $p$ is an odd prime and $n$
is an integer greater than one. \noindent Since
$$\ord_p(m!)=\sum_{i=1}^{\infty}\left\lfloor\frac m{p^i}\right\rfloor
\qquad\t{for every}\ m=0,1,2,\ldots$$ (a well-known fact in number
theory), the inequality in Theorem \ref{mainthm} can be restated
as
$$\exp_p({\rm SU}(n))\ge n-1+\sum_{i=2}^{\infty}\l\lfloor\frac n{p^i}\r\rfloor,$$
a nice improvement of (\ref{largethm}).

\section{Outline of proof}\label{outlinesec}
In this section we present the deduction of Theorem \ref{mainthm}
from Theorem \ref{combthm}, which will then be proved in Section
\ref{numthsec}. We also present some comments regarding the extent
to which Theorem \ref{mainthm} is sharp.

Let $p$ be any prime. In \cite{DM}, the first author and Mahowald
defined the ($p$-primary) $v_1$-periodic homotopy groups
$\vp_*(X;p)$ of a topological space $X$ and proved that if $X$ is
a sphere or compact Lie group, such as ${\rm SU}(n)$, each group
$\vp_i(X;p)$ is a direct summand of some actual homotopy group
$\pi_j(X)$. See also \cite{sur} for another expository account of
$v_1$-periodic homotopy theory.

In \cite[1.4]{DSU} and \cite[1.1a]{BDSU}, it was proved that if
$p$ is odd, or if $p=2$ and $n$ is odd, then there is an
isomorphism
\begin{equation}\label{vpSUn}\vp_{2k}
({\rm SU}(n);p)\cong\bz/p^{e_p(n,k)}\Z
\end{equation}
for all $k\ge n$, where $e_p(n,k)$ is as defined in \ref{edef} and
we use $\bz/m\Z$ to denote the additive group of residue classes
modulo $m$. Thus, unless $p=2$ and $n$ is even, for any integer
$k\ge n$, we have
$$\exp_p({\rm SU}(n))\ge e_p(n,k),$$
establishing Proposition \ref{redn} in these cases. The situation
when $p=2$ and $n$ is even is somewhat more technical, and will be
discussed in Section \ref{p=2sec}.

Next we show that Theorem \ref{combthm} implies Theorem
\ref{numththm}.
\begin{pf*}{Proof of Theorem \ref{numththm}}
(i) By a well-known property of Stirling numbers of the second
kind (cf. \cite [pp.\,125-126]{LW}),
$$m!S(kh(p-1)p^{\al}+n-1,m)=\sum_{j=0}^m\bi mj(-1)^{m-j}j^{kh(p-1)p^{\al}+n-1}$$
for any $k\in\N$. Thus
$$(-1)^m m!\sum_{k=0}^l\bi lk(-1)^kS(kh(p-1)p^{\al}+n-1,m)=\Sigma_1+\Sigma_2,$$
where
\begin{eqnarray*}\Sigma_1&=&
\sum_{k=0}^l\bi lk(-1)^kp^{n-1+kh(p-1)p^{\al}} \sum_{j\equiv
0\,(\mo\ p)}\bi mj(-1)^j\l(\f jp\r)^{n-1+kh(p-1)p^{\al}}
\end{eqnarray*}
and
\begin{eqnarray*}\Sigma_2&=&
\sum_{j\not\eq0\,(\mo\ p)}\bi mj(-1)^j\sum_{k=0}^l\bi
lk(-1)^kj^{n-1+kh(p-1)p^{\al}}
\\&=&\sum_{j\not\eq0\,(\mo\ p)}\bi mj(-1)^jj^{n-1}\l(1-j^{h(p-1)p^{\al}}\r)^l.
\end{eqnarray*}
Clearly $\ord_p(\Sigma_1)\ge n-1+\ord_p(\lfloor m/p\rfloor!)$ by
Theorem \ref{combthm}, and $\ord_p(\Sigma_2)\ge l(\al+1)$ by
Euler's theorem in number theory. Therefore the first part of
Theorem \ref{numththm} holds.

(ii) Observe that $$N+1-(n-1)>\f{n}{p(p-1)}=\sum_{i=2}^{\infty}\f
n{p^i}
>\sum_{i=2}^{\infty}\l\lfloor\f n{p^i}\r\rfloor=\ord_p\l(\l\lfloor\f np\r\rfloor!\r).$$

By part (i) in the case $l=h=1$ and $\al=L\gs N$, if $m\ge n$ then
\begin{eqnarray*}&&\ord_p\l(m!S(n-1,m)-m!S((p-1)p^{L}+n-1,\,m)\r)
\\&\ &\qquad\ge n-1+\ord_p\l(\l\lfloor\f np\r\rfloor!\r).
\end{eqnarray*}
Since $S(n-1,m)=0$ for $m\ge n$, we finally have
$$e_p(n,(p-1)p^{L}+n-1)\ge n-1+\ord_p\l(\l\lfloor\f np\r\rfloor!\r)$$
as required.\end{pf*}

\smallskip
The following proposition, although not needed for our main
results, sheds more light on the large exponents $N$ and $L$ which
appear in Theorem \ref{numththm}(ii), and is useful in our
subsequent exposition.

\begin{prop}\label{indptprop} Let $p$ be a prime and $n>1$ be an integer.
Then there exists an integer $N_0\ge0$, effectively computable in
terms of $p$ and $n$, such that $e_p(n,(p-1)p^L+n-1)$ has the same
value for all $L\ge N_0$.
\end{prop}

\begin{pf} For integers $m\gs n$ and $L\ge0$, we write
$$(-1)^mm!S((p-1)p^L+n-1,m)=\sum_{j=0}^m\bi mj(-1)^jj^{(p-1)p^L+n-1}=S_m+S'_{m,L}+S''_{m,L},$$
where
\begin{eqnarray*}S_m&=&\sum_{j\not\equiv0\,(\mo\ p)}\binom mj(-1)^jj^{n-1},\\
S'_{m,L}&=&\sum_{j\not\equiv0\,(\mo\ p)}\binom mj (-1)^jj^{n-1}(j^{(p-1)p^L}-1),\\
S''_{m,L}&=&\sum_{j\equiv0\,(\mo\ p)}\binom mj(-1)^j
j^{(p-1)p^L+n-1}.
\end{eqnarray*}
Note that both $S'_{m,L}$ and $S''_{m,L}$ are divisible by
$p^{L+1}$.

Assume that $S_n,S_{n+1},\ldots$ are not all zero. (This will be
shown later.) Then $L_0=\min_{m\ge n}\ord_p(S_m)$ is finite. Let
$m_0\gs n$ satisfy $\ord_p(S_{m_0})=L_0$. Whenever $L\ge L_0$, we
have $\ord_p(S_m+S'_{m,L}+S''_{m,L})\ge L_0$ for every $m\ge n$,
and equality is attained for $m=m_0$. Thus, if $L\gs L_0$ then
\begin{eqnarray*}e_p(n,(p-1)p^L+n-1)&=&\min_{m\ge n}\ord_p(m!S((p-1)p^L+n-1,m))
\\&=&\min_{m\ge n}\ord_p(S_m+S'_{m,L}+S''_{m,L})=L_0.
\end{eqnarray*}
Although $L_0$ is finite, it may not be effectively computable.
Instead of $L_0$ we use the $p$-adic order $N_0$ of the first
nonzero term in the sequence $S_n,S_{n+1},\ldots$. This $N_0$ is
computable, also $e_p(n,(p-1)p^L+n-1)=L_0$ for all $L\ge N_0$
since $N_0\ge L_0$.

To complete the proof, we must show that $S_m$ is nonzero for some
$m\ge n$. First note that this is clearly true for $p=2$ since
then $S_m$ is a sum of negative terms. If $p$ is odd and $S_m=0$
for all $m\ge n$, then $e_p(n,(p-1)p^L+n-1)=\min_{m\gs
n}\ord_p(S'_{m,L}+S''_{m,L})\ge L+1$ for any $L\ge0$. By
(\ref{vpSUn}), this would imply that $\vp_*({\rm SU}(n);p)$ has
elements of arbitrarily large $p$-exponent. However, this is not
true, for in \cite[5.8]{DSU}, it was shown that the $v_1$-periodic
$p$-exponent of ${\rm SU}(n)$ does not exceed
$e:=\lfloor(n-1)(1+(p-1)^{-1}+(p-1)^{-2})\rfloor$; i.e., for this
$e$,  $p^e\vp_*({\rm SU}(n);p)=0$.
\end{pf}

In the remainder of this section and in Section \ref{p=2sec}, once
a prime $p$ and an integer $n>1$ is given, $L$ will refer to any
integer not smaller than $\max\{N,N_0\}$ where $N$ and $N_0$ are
described in Theorem \ref{numththm}(ii) and the proof of
Proposition \ref{indptprop} respectively.

\smallskip

We now comment on the extent to which Theorem \ref{mainthm} might
be sharp. In Table \ref{gptbl}, we present, for $p=3$ and a
representative set of values of $n$, three numbers. The first,
labelled $\exp_3(v_1^{-1}{\rm SU}(n))$, is the largest value of
$e_3(n,k)$ over all values of $k\ge n$; thus it is the largest
exponent of the 3-primary $v_1$-periodic homotopy groups of ${\rm
SU}(n)$. The second number in the table is the exponent of the
$v_1$-periodic homotopy group on which we have been focusing,
which, at least in the range of this table, is equal to or just
slightly less than the maximal exponent. The third number is the
nice estimate for this exponent given by Theorem
\ref{numththm}(ii).

\begin{table}[h]
\begin{center}
\caption{Comparison of exponents when $p=3$}\label{gptbl}
\begin{tabular}{c|ccc}
$n$&$\exp_3(v_1^{-1}{\rm SU}(n))$&$e_3(n,2\cdot3^L+n-1)$&$n-1+\ord_3(\lfloor n/3\rfloor!)$\\
\hline
19&21&20&20\\
20&22&21&21\\
21&22&22&22\\
22&25&25&23\\
23&26&26&24\\
24&28&28&25\\
25&29&28&26\\
26&30&30&27\\
27&31&31&30\\
28&32&32&31\\
29&34&32&32\\
30&34&33&33\\
31&34&34&34\\
32&35&35&35\\
33&37&37&36\\
34&38&37&37\\
35&39&39&38\\
36&41&41&40\\
37&42&41&41\\
38&43&42&42\\
39&43&43&43\\
40&45&44&44\\
41&45&45&45
\end{tabular}
\end{center}
\end{table}

Note that, for more than half of the values of $n$ in the table,
the largest group $\vp_{2k}({\rm SU}(n);3)$ occurs when
$k=2\cdot3^L+n-1$. In the worst case in the table, $n=29$,
detailed {\tt Maple} calculations suggest that if $k\ge 29$ and
$k\equiv10\ (\mo \ 18)$, then $$e_3(29,k)=\min\{\ord_3(k-28-8\cdot
3^{20})+12,\,34\}.$$ Shifts (as by $8\cdot3^{20}$) were already
noted in \cite[p.\,543]{DSU}. Note also that for more than half of
the cases in the table, our estimate for $e_3(n,2\cdot3^L+n-1)$ is
sharp, and it never misses by more than 3.

 The big question for
topologists, though, is whether the $v_1$-periodic $p$-exponent
agrees (or almost agrees) with the actual homotopy $p$-exponent.
The fact that they agree for $S^{2n+1}$ when $p$ is an odd prime
(\cite{CMN1}, \cite{Gray}) leads the first author to conjecture
that they also agree for ${\rm SU}(n)$ if $p\not=2$, but we have
no idea how to prove this. Theriault (\cite{Th1}, \cite{Th2}) has
made good progress in proving that some of the first author's
lower bounds for $p$-exponents of certain exceptional Lie groups
are sharp.

\section{Proof of Theorem \ref{combthm}}\label{numthsec}
In this section, we prove Theorem \ref{combthm}, which we have
already shown to imply Theorem \ref{mainthm}.

\begin{lem}\label{Lem3.1}  Let $p$ be any prime, and let $\al,n\in\N$ and $r\in\Z$.
Then
$$\ord_p\(\sum_{k\eq r\,(\mo\ p^{\al})}\bi nk(-1)^k\)
\gs\ord_p\l(\l\lfloor\f n{p^{\al-1}}\r\rfloor!\r) =\l\lfloor\f
n{p^{\al}}\r\rfloor+\ord_p\l(\l\lfloor\f
n{p^{\al}}\r\rfloor!\r).$$
\end{lem}

\begin{pf} The equality is easy, for,
\begin{eqnarray*}
\ord_p\l(\l\lfloor\f n{p^{\al-1}}\r\rfloor!\r)
&=&\sum_{i=1}^{\infty}\l\lfloor\f{\lfloor
n/p^{\al-1}\rfloor}{p^i}\r\rfloor
=\sum_{j=\al}^{\infty}\l\lfloor\f n{p^j}\r\rfloor
\\&=&\l\lfloor\f n{p^{\al}}\r\rfloor
+\sum_{i=1}^{\infty}\l\lfloor\f{\lfloor
n/p^{\al}\rfloor}{p^i}\r\rfloor =\l\lfloor\f
n{p^{\al}}\r\rfloor+\ord_p\l(\l\lfloor\f n{p^{\al}}\r\rfloor!\r).
\end{eqnarray*}

When $\al=0$ or $n<p^{\al-1}$, the desired inequality is obvious.

Now let $\al>0$ and $m=\lfloor n/p^{\al-1}\rfloor\gs1$. Observe
that
$$\ord_p(m!)=\sum_{i=1}^{\infty}\l\lfloor\f m{p^i}\r\rfloor
<\sum_{i=1}^{\infty}\f m{p^i}=\f mp\sum_{j=0}^{\infty}\f1{p^j} =\f
mp\cdot\f1{1-p^{-1}}=\f{m}{p-1}.$$ Thus $(p-1)\ord_p(m!)\ls m-1$,
and hence
$$\ord_p(m!)\ls\l\lfloor\f{m-1}{p-1}\r\rfloor
=\l\lfloor \f{n/p^{\al-1}-1}{p-1}\r\rfloor
=\l\lfloor\f{n-p^{\al-1}}{\varphi(p^{\al})}\r\rfloor,$$ where
$\varphi$ is Euler's totient function. By a result of Weisman
\cite{W77},
$$\ord_p\(\sum_{k\eq r\,(\mo\ p^{\al})}\bi nk(-1)^k\)
\gs\l\lfloor\f{n-p^{\al-1}}{\varphi(p^{\al})}\r\rfloor.$$
(Weisman's proof is complicated, but an easy induction proof
appeared in \cite{S05}.) So we have the desired inequality.
\end{pf}
Now we restate Lemma 2.1 of Sun \cite{S05}, which will be used
later.
\begin{lem} $(\cite{S05})$ Let $m$ and $n$ be positive integers, and let $f(x)$ be a function from
$\Z$ to a field. Then, for any $r\in\Z$, we have
$$\sum_{k=0}^n\binom nk (-1)^kf\l(\l\lfloor\frac{k-r}m\r\rfloor\r)=\sum_{k\equiv \rbar(\mo\ m)}\binom{n-1}k
(-1)^{k-1}\Delta f\l(\frac{k-\rbar}m\r),$$ where $\rbar=r-1+m$ and
$\Delta f(x)=f(x+1)-f(x)$.\label{Sunlem}
\end{lem}
\begin{lem}\label{Lem3.3}  Let $m,n\in\Z^+$ and $r\in\Z$, and let $f(x)$
be a complex-valued function defined on $\Z$. Then we have
\begin{eqnarray*}&&\sum_{k\eq r\,(\mo\ m)}\bi nk(-1)^kf\l(\f{k-r}m\r)
-f\l(\l\lfloor\f{n-r}m\r\rfloor\r)\sum_{k\eq r\,(\mo\ m)}\bi
nk(-1)^k
\\&\ &=-\sum_{j=0}^{n-1}\bi nj\sum_{i\eq r\,(\mo\ m)}\bi ji(-1)^i
\sum_{k\eq r_j\,(\mo\ m)}\bi{n-j-1}{k}(-1)^{k} \Delta
f\l(\f{k-r_j}m\r),
\end{eqnarray*}
where $r_j=r-j+m-1$.
\end{lem}

\begin{pf} Let $\zeta$ be a primitive $m$th root of unity.
Clearly
$$\sum_{s=0}^{m-1}\zeta^{(k-r)s}=\cases \sum_{s=0}^{m-1}1=m&\t{if}\ k\eq r\,(\mo\ m),\\
\f{1-\zeta^{(k-r)m}}{1-\zeta^{k-r}}=0&\t{otherwise}.\endcases$$
Thus
\begin{eqnarray*}&&\sum_{k\eq r\,(\mo\ m)}\bi nk(-1)^kf\l(\f{k-r}{m}\r)
\\&=&\sum_{k=0}^n\(\f1{m}\sum_{s=0}^{m-1}\zeta^{(k-r)s}\)
\bi nk(-1)^kf\(\l\lfloor\f{k-r}m\r\rfloor\)
=\f1m\sum_{s=0}^{m-1}\zeta^{-rs}c_s,
\end{eqnarray*}
where
$$c_s=\sum_{k=0}^n\bi nk(-\zeta^s)^kf\(\l\lfloor\f{k-r}m\r\rfloor\).$$

Observe that
\begin{eqnarray*}c_s&=&\sum_{k=0}^n\bi nk((1-\zeta^s)-1)^kf\(\l\lfloor\f{k-r}m\r\rfloor\)
\\&=&\sum_{k=0}^n\bi nk\sum_{j=0}^k\bi kj(1-\zeta^s)^j(-1)^{k-j}f\(\l\lfloor\f{k-r}m\r\rfloor\)
\\&=&\sum_{j=0}^n\bi nj(1-\zeta^s)^j\sum_{k=j}^n\bi{n-j}{k-j}(-1)^{k-j}
f\(\l\lfloor\f{k-r}m\r\rfloor\)
\\&=&\sum_{j=0}^n\bi nj(1-\zeta^s)^j\sum_{k=0}^{n-j}\bi{n-j}{k}(-1)^{k}
f\(\l\lfloor\f{k-(r-j)}m\r\rfloor\).
\end{eqnarray*}
Applying Lemma \ref{Sunlem}, we find that
\begin{eqnarray*}&&c_s-(1-\zeta^s)^nf\l(\l\lfloor\f{n-r}m\r\rfloor\r)
\\&=&\sum_{j=0}^{n-1}
\bi nj(1-\zeta^s)^j\sum_{k\eq r_j\, (\mo\
m)}\bi{n-j-1}{k}(-1)^{k-1} \Delta f\l(\f{k-r_j}m\r).
\end{eqnarray*}

In view of the above, it suffices to note that
$$\sum_{s=0}^{m-1}\f{\zeta^{-rs}}m(1-\zeta^s)^j=
\sum_{i=0}^j\bi ji\f{(-1)^i}m\sum_{s=0}^{m-1}\zeta^{s(i-r)}
=\sum_{i\eq r\,(\mo\ m)}\bi ji(-1)^i.$$ This concludes the proof.
\end{pf}

With help of Lemmas \ref{Lem3.1} and \ref{Lem3.3}, we are able to
prove the following equivalent version of Theorem \ref{combthm}.

\begin{thm}\label{combthm'} Let $p$ be a prime, and let
$\al,l,n\in\N$. Then for any $r\in\Z$ we have
$$\ord_p\(\sum_{k\eq r\,(\mo\ p^{\al})}\bi nk(-1)^k\l(\f{k-r}{p^{\al}}\r)^l\)
\gs\ord_p\l(\l\lfloor\f n{p^{\al}}\r\rfloor!\r).$$
\end{thm}

\begin{pf}
We use induction on $l$.

In the case $l=0$, the desired result follows from Lemma 3.1.

 Now let $l>0$ and assume the result for smaller values of $l$.
We use induction on $n$ to prove the inequality in Theorem
\ref{combthm'}.

The case $n=0$ is trivial. So we now let $n>0$ and assume that the
inequality holds with smaller values of $n$. Observe that
\begin{eqnarray*}&&\sum_{k\eq r\,(\mo\ p^{\al})}\bi nk(-1)^k\l(\f{k-r}{p^{\al}}\r)^l
\\&=&\sum_{k\eq r\,(\mo\ p^{\al})}\(\bi{n-1}k+\bi{n-1}{k-1}\)(-1)^k\l(\f{k-r}{p^{\al}}\r)^l
\\&=&\sum_{k\eq r\,(\mo\ p^{\al})}\bi {n-1}k(-1)^k\l(\f{k-r}{p^{\al}}\r)^l
\\&\ &-\sum_{k'\eq r-1\,(\mo\ p^{\al})}\bi{n-1}{k'}(-1)^{k'}\l(\f{k'-(r-1)}{p^{\al}}\r)^l.
\end{eqnarray*}
In view of this, if $p^{\al}$ does not divide $n$, then, by the
induction hypothesis for $n-1$, we have
\begin{eqnarray*}&&\ord_p\(\sum_{k\eq r\,(\mo\ p^{\al})}\bi nk(-1)^k\l(\f{k-r}{p^{\al}}\r)^l\)
\\&&\ \ \gs\ord_p\l(\l\lfloor\f{n-1}{p^{\al}}\r\rfloor\r)
=\ord_p\l(\l\lfloor\f{n}{p^{\al}}\r\rfloor\r).
\end{eqnarray*}

\medskip

Below we let $p^{\al}\mid n$ and set $m=n/p^{\al}$.

{\it Case 1}. $r\eq 0\ (\mo\ p^{\al})$. In this case,
\begin{eqnarray*}&&\f1{m!}\sum_{k\eq r\,(\mo\ p^{\al})}
\bi nk(-1)^k\l(\f k{p^{\al}}-\f r{p^{\al}}\r)
\l(\f{k-r}{p^{\al}}\r)^{l-1}
\\&=&\f{n/p^{\al}}{m!}\sum_{k\eq r\,(\mo\ p^{\al})}
\bi{n-1}{k-1}(-1)^k\l(\f{k-r}{p^{\al}}\r)^{l-1}
\\&\ &-\f{r/p^{\al}}{m!}\sum_{k\eq r\,(\mo\ p^{\al})}
\bi{n}{k}(-1)^k\l(\f{k-r}{p^{\al}}\r)^{l-1}
\\&=&\f1{\lfloor (n-1)/p^{\al}\rfloor!}\sum_{k\eq r-1\,(\mo\ p^{\al})}\bi{n-1}k(-1)^{k+1}
\l(\f{k-(r-1)}{p^{\al}}\r)^{l-1}
\\&\ &-\f{r/p^{\al}}{\lfloor n/p^{\al}\rfloor!}
\sum_{k\eq r\,(\mo\
p^{\al})}\bi{n}{k}(-1)^k\l(\f{k-r}{p^{\al}}\r)^{l-1}.
\end{eqnarray*}
Thus, by the induction hypothesis for $l-1$,
$$\f1{m!}\sum_{k\eq r\,(\mo\ p^{\al})}\bi nk(-1)^k\l(\f {k-r}{p^{\al}}\r)^l$$
is a $p$-integer (i.e., its denominator is relatively prime to
$p$) and hence the desired inequality follows.

{\it Case 2}. $r\not\eq 0\ (\mo\ p^{\al})$. Note that $\sum_{i\eq
r\,(\mo\ p^{\al})}\bi 0i(-1)^i=0$. Also,
$$\ord_p\(\sum_{k\eq r\,(\mo\ p^{\al})}\bi nk(-1)^k\)
\ge\ord_p\l(\f n{p^{\al-1}}!\r) =m+\ord_p(m!)$$ by Lemma
\ref{Lem3.1}. Thus, in view of Lemma \ref{Lem3.3},
 it suffices to show that if $0<j<n$
then the $p$-adic order of
$$\sigma_j=\bi nj\sum_{i\eq r\,(\mo\ p^{\al})}\bi ji(-1)^i
\sum_{k\eq r_j\,(\mo\ p^{\al})}\bi{n-j-1}{k}(-1)^{k} \Delta
f\l(\f{k-r_j}{p^{\al}}\r)$$ is at least $\ord_p(m!)$, where
$r_j=r-j+p^{\al}-1$ and $f(x)=x^l$.

Let $0<j\ls n-1$ and write $j=p^{\al}s+t$, where $s,t\in\N$ and
$t<p^{\al}$. Note that
$$\l\lfloor\f{j}{p^{\al}}\r\rfloor=s\quad\ \t{and}\quad\
\l\lfloor\f{n-j-1}{p^{\al}}\r\rfloor=\l\lfloor
m-s-\f{t+1}{p^{\al}}\r\rfloor=m-s-1.$$ Since $\Delta
f(x)=(x+1)^l-x^l=\sum_{i=0}^{l-1}\bi lix^i$, by Lemma 3.1 and the
induction hypothesis with respect to $l$, we have
\begin{eqnarray*}\ord_p(\sigma_j)&=&\ord_p\bi nj
+\ord_p\(\sum_{i\eq r\,(\mo\ p^{\al})}\bi ji(-1)^i\)
\\&\ &+\ord_p\(\sum_{k\eq r_j\,(\mo\ p^{\al})}\bi{n-j-1}{k}(-1)^{k}
\Delta f\l(\f{k-r_j}{p^{\al}}\r)\)
\\&\gs&\ord_p\bi nj+(s+\ord_p(s!))
+\ord_p((m-s-1)!)
\\&=&\ord_p\bi nj+s+\ord_p(s!)-\ord_p\(\prod_{i=0}^s(m-i)\)
+\ord_p(m!)
\\&=&\ord_p\bi {p^{\al}m}{p^{\al}s+t}-\ord_p\bi ms+s-\ord_p(m-s)
+\ord_p(m!).
\end{eqnarray*}
When $t=0$ (i.e., $j=p^{\al}s$) we have the stronger inequality
$$\ord_p(\sigma_j)\ge\ord_p\bi {p^{\al}m}{p^{\al}s}-\ord_p\bi ms+s
+\ord_p(m!),$$ because
\begin{eqnarray*}&&\ord_p\(\sum_{k\eq r_j\,(\mo\ p^{\al})}\bi{n-j-1}{k}(-1)^{k} \Delta
f\l(\f{k-r_j}{p^{\al}}\r)\)
\\&=&\ord_p\(\sum_{k=0}^{n-j}\bi{n-j}k(-1)^kf\(\l\lfloor\f{k-(r-j)}{p^{\al}}\r\rfloor\)\)
\ \ (\t{by Lemma \ref{Sunlem}})
\\&=&\ord_p\(\sum_{i=0}^{p^{\al}-1}
\sum_{k-(r-j)\eq i\, (\mo\
p^{\al})}\bi{n-j}k(-1)^k\(\f{k-(r-j)-i}{p^{\al}}\)^l\)
\\&\ge&\ord_p\(\f{n-j}{p^{\al}}!\)=\ord_p((m-s)!)
\\& &\qquad\qquad\quad\ (\t{by the induction hypothesis with respect to}\ n).
\end{eqnarray*}

Observe that
\begin{eqnarray*}\ord_p\bi{p^{\al}m}{p^{\al}s}
&=&\sum_{i=1}^{\al}\(\f{p^{\al}m}{p^i}-\f{p^{\al}s}{p^i}-\f{p^{\al}(m-s)}{p^i}\)
\\&\ &+\sum_{i=\al+1}^{\infty}\(\l\lfloor\f{p^{\al}m}{p^i}\r\rfloor
-\l\lfloor\f{p^{\al}s}{p^i}\r\rfloor-\l\lfloor\f{p^{\al}(m-s)}{p^i}\r\rfloor\)
\\&=&\sum_{i=1}^{\infty}\(\l\lfloor\f{m}{p^i}\r\rfloor
-\l\lfloor\f{s}{p^i}\r\rfloor-\l\lfloor\f{m-s}{p^i}\r\rfloor\)=\ord_p\bi
ms.
\end{eqnarray*}
Thus, when $t=0$ we have $\ord_p(\sigma_j)\gs
s+\ord_p(m!)\gs\ord_p(m!)$.

Define $\ord_p(a/b)=\ord_p(a)-\ord_p(b)$ if $a,b\in\Z$ and $a$ is
not divisible by $b$. If $t>0$ then
\begin{eqnarray*}&&\ord_p\bi{p^{\al}m}{p^{\al}s+t}-\ord_p\bi ms
\\&=&\ord_p\f{\bi{p^{\al}m}{p^{\al}s+t}}{\bi{p^{\al}m}{p^{\al}s}}
=\ord_p\f{(p^{\al}s)!(p^{\al}(m-s))!}{(p^{\al}s+t)!(p^{\al}(m-s)-t)!}
\\&=&\ord_p\f{p^{\al}(m-s)}{p^{\al}s+t}
+\ord_p\prod_{0<i<t}\f{p^{\al}(m-s)-i}{p^{\al}s+i}.
\end{eqnarray*}
For $0<i<p^{\al}$, clearly
$$\ord_p(p^{\al}(m-s)-i)=\ord_p(p^{\al}s+i)=\ord_p(i)<\al.$$
Therefore, when $0<t<p^{\al}$ we have
$$\ord_p\bi{p^{\al}m}{p^{\al}s+t}-\ord_p\bi ms
=\ord_p(m-s)+\al-\ord_p(p^{\al}s+t)>\ord_p(m-s)$$ and hence
$\ord_p(\sigma_j)>s+\ord_p(m!)\gs \ord_p(m!)$. This concludes the
analysis of the second case.

 The proof of Theorem \ref{combthm'} is now complete.
\end{pf}

 Note that, in the proof of Theorem \ref{combthm'},
the technique used to handle the first case is of no use in the
second case, and vice versa. Thus, the distinction of the two
cases is important.

\begin{cor}\label{DS}  Let $p$ be a prime, and let $\al,l,n\in\N$ and $r\in\Z$.
Then we have
$$\ord_p\(\sum_{k\eq r\,(\mo\ p^{\al})}\bi nk(-1)^k\bi{(k-r)/p^{\al}}l\)
\gs\ord_p\l(\l\lfloor\f n{p^{\al}}\r\rfloor!\r)-\ord_p(l!).$$
\end{cor}
\begin{pf} Simply apply Theorem \ref{combthm} with $f(x)=l!\bi xl\in\Z[x]$.
\end{pf}

\section{Changes when $p=2$ and $n$ is even}\label{p=2sec}
When $p=2$ and $n$ is even, the relationship between
$\vp_{2k}({\rm SU}(n);p)$ and $e_p(n,k)$ (with $k\ge n$) is not so
simple as in (\ref{vpSUn}). As described in \cite{BDSU} and
\cite{DP}, there is a spectral sequence converging to $\vp_*({\rm
SU}(n);p)$ and satisfying $E_2^{1,\,2k+1}({\rm SU}(n))\cong
\bz/p^{e_p(n,k)}\Z$. If $p$ or $n$ is odd, the spectral sequence
necessarily collapses and $\vp_{2k}({\rm SU}(n);p)\cong
E_2^{1,\,2k+1}$. (Here we begin abbreviating $E_r^{*,\,*}({\rm
SU}(n))$ just as $E_r^{*,\,*}$.) If $p=2$ and $n$ is even, there
are two ways in which the corresponding summand of $\vp_{2k}({\rm
SU}(n);2)$ may differ from this.

It is conceivable that there could be an extension in the spectral
sequence, which would make the exponent of the homotopy group 1
larger than that of $E_\infty^{1,\,2k+1}$.  However, as observed
in \cite[6.2(1)]{DP}, it is easily seen that this does not happen.

It is also conceivable that the differential
$d_3:E_3^{1,\,2k+1}\to E_3^{4,\,2k+3}\cong\bz/2\Z\oplus\bz/2\Z$
could be nonzero, which would make the exponent of $\vp_{2k}({\rm
SU}(n);2)$ equal to $e_2(n,k)-1$. This is the reason for the $-1$
at the end of Proposition \ref{redn}. By \cite[1.6]{DP}, if
$n\equiv0\ (\mo\ 4)$ and $k=2^L+n-1$, then $d_3: E_3^{1,\,2k+1}\to
E_3^{4,\,2k+3}$ must be 0.

Now suppose $n\equiv2\ (\mo\ 4)$. If $n=2$, then
$n-1+\ord_2(\lfloor n/2\rfloor!)=1<\exp_2({\rm SU}(n))$ since
$\pi_6({\rm SU}(2))\cong\Z/12\Z$ (cf. \cite{To}). Below we let
$n>2$, hence $n/2+1$ is even and not larger than $n-1$. As first
noted in \cite[1.1]{BDSU} and restated in \cite[6.5]{DP}, for
$k=2^L+n-1$, the differential $d_3:E_3^{1,\,2k+1}\to
E_3^{4,\,2k+3}$ is nonzero if and only if
$$e_2(n,2^L+n-1)=e_2(n-1,2^L+n-1)+n-1.$$
We show at the end of the section that
\begin{equation}\label{factorial}e_2(n-1,2^L+n-1)=\ord_2((n-1)!).\end{equation}
Thus, if the above $d_3$ is nonzero, then
$e_2(n,2^L+n-1)=n-1+\ord_2((n-1)!)$ and hence
$$\exp_2({\rm SU}(n))\ge e_2(n,2^L+n-1)-1= n-1+\ord_2((n-1)!)-1
\ge n-1+\ord_2(\lfloor n/2\rfloor!),$$ as claimed in Theorem
\ref{mainthm}.

\begin{pf*}{Proof of $(\ref{factorial})$}
Putting $p=2$, $\al=L$, $l=h=1$ and $m=n-1$ in the first part of
Theorem \ref{numththm}, we get that
\begin{eqnarray*}&&\ord_2\l((n-1)!S(n-1,n-1)-(n-1)!S(2^{L}+n-1,n-1)\r)
\\&&\quad\ge n-1+\ord_2\l(\l\lfloor\f{n-1}2\r\rfloor!\r)\gs n-1>\ord_2((n-1)!).
\end{eqnarray*}
Therefore $\ord_2((n-1)!S(2^{L}+n-1,n-1))=\ord_2((n-1)!)$. On the
other hand, by the second part of Theorem \ref{numththm},
$\ord_2(m!S(2^{L}+n-1,m)) \ge n-1+\ord_2(\lfloor n/2\rfloor!)$ for
all $m\ge n$. So we have (\ref{factorial}).
\end{pf*}

\section{Strengthening and sharpness of Theorem \ref{combthm'}}\label{strongsec}
In this section, we give an example illustrating the extent to
which Theorem \ref{combthm'} is sharp when $r=0$, which is the
situation that is used in our application to topology. Then we
show in Theorem \ref{strthm} that the lower bound in Theorem
\ref{combthm'} can sometimes be increased slightly.

We begin with a typical example of Theorem \ref{combthm'}. Let
$p=\al=2$, $r=0$ and $n=100$. Then $\lfloor n/p^{\al}\rfloor=25$
and $\ord_p(\lfloor n/p^{\al}\rfloor!)=22$. For $l\ge25$, set
\begin{equation*}\label{fdef}\delta(l)=\ord_2\(\sum_{k\equiv 0\,(\mo\ 4)}
\binom nk\l(\frac{k}{4}\right)^{l}\)-22.\end{equation*} The range
$l\ge\lfloor n/p^\a\rfloor=25$ is that in which we feel Theorem
\ref{combthm'} to be very strong. (See Remark \ref{rmks}(2).)
Clearly $\delta(l)$ measures the amount by which the actual
$p$-adic order of the sum in Theorem \ref{combthm'} exceeds our
bound for it. The values of $\delta(l)$ for $25\le l\le 45$ are
given in order as
$$0,0,0,0,2,3,2,4,1,1,1,1,2,2,4,1,0,0,0,0,3.$$

When $r=0$ and in many other situations, Theorem \ref{combthm'}
appears to be sharp for infinitely many values of $l$.

Before presenting our strengthening of Theorem \ref{combthm'} we
need some notation. For $a\in\Z$ and $m\in\Z^+$, we let $\{a\}_m$
denote the least nonnegative residue of $a$ modulo $m$. Given a
prime $p$, for any $a,b\in\N$ we let $\tau_p(a,b)$ represent the
number of carries occurring in the addition of $a$ and $b$ in base
$p$; actually
$$\tau_p(a,b)=\sum_{i=1}^{\infty}\(\l\lfloor\f {a+b}{p^i}\r\rfloor
-\l\lfloor\f {a}{p^i}\r\rfloor-\l\lfloor\f {b}{p^i}\r\rfloor\)
=\ord_p\bi{a+b}a$$ as observed by E. Kummer.

Here is our strengthening of Theorem \ref{combthm'}. The right
hand side is the amount by which the bound in Theorem
\ref{combthm'} can be improved. This amount does not exceed $\al$,
by the definition of $\tau_p$. In Table \ref{comtbl}, we
illustrate this amount when $p=3$ and $\al=2$.

\begin{thm}\label{strthm} Let $p$ be a prime, and let $\al,l,n\in\N$.
Then, for all $r\in\Z$, we have
\begin{eqnarray*}&&\ord_p\(\sum_{k\equiv r\,(\mo\ p^\a)}
\binom nk(-1)^k\left(\frac{k-r}{p^\a}\right)^l\)
-\ord_p\left(\l\lfloor\frac n{p^\a}\r\rfloor!\right)
\\&&\quad\ge \tau_p(\{r\}_{p^{\al}},\{n-r\}_{p^{\al}})
=\ord_p\binom {\{r\}_{p^{\al}}+\{n-r\}_{p^{\al}}}
{\{r\}_{p^{\al}}}.
\end{eqnarray*}
\end{thm}

\begin{pf}
We use induction on $n$.

In the case $n=0$, whether $r\eq0\ (\mo\ p^{\al})$ or not, the
desired result holds trivially.

Now let $n>0$ and assume the corresponding result for $n-1$.
Suppose that $\tau_p(\{r\}_{p^{\al}},\{n-r\}_{p^{\al}})>0$. Then
neither $r$ nor $n-r$ is divisible by $p^{\al}$.

Set
$$R=\f1{\lfloor n/p^{\al}\rfloor!}
\sum_{k\eq r\,(\mo\ p^{\al})}\bi
nk(-1)^k\l(\f{k-r}{p^{\al}}\r)^l$$ and
$$R'=\f{n/p^{\al}}{\lfloor n/p^{\al}\rfloor!}\label{R'eq}
\sum_{k\eq r-1\,(\mo\ p^{\al})}\bi
{n-1}k(-1)^k\l(\f{k-(r-1)}{p^{\al}}\r)^{l}.$$ Clearly
\begin{eqnarray*}R'&=&-\f{n/p^{\al}}{\lfloor n/p^{\al}\rfloor!}\nonumber
\sum_{k\eq r\,(\mo\ p^{\al})}\bi
{n-1}{k-1}(-1)^k\l(\f{k-r}{p^{\al}}\r)^{l}
\\&=&-\f1{\lfloor n/p^{\al}\rfloor!}
\sum_{k\eq r\,(\mo\ p^{\al})}\bi {n}k(-1)^k\f
k{p^{\al}}\l(\f{k-r}{p^{\al}}\r)^{l},
\end{eqnarray*}
and thus
\begin{equation*}\f{r}{p^{\al}}R+R'=-\f{1}{\lfloor n/p^{\al}\rfloor!}
\sum_{k\eq r\,(\mo\ p^{\al})}\bi
nk(-1)^k\l(\f{k-r}{p^{\al}}\r)^{l+1}.
\end{equation*}
This is a $p$-integer by Theorem \ref{combthm'}; therefore
$\ord_p(rR+p^{\al}R')\ge\al$.

Let $\b=\ord_p(n)$. We consider three cases.

{\it Case 1}. $\beta\ge\al$. In this case, $\lfloor
n/p^{\al}\rfloor!/(n/p^{\al})=\lfloor(n-1)/p^{\al}\rfloor!$ and
hence $R'$ is a $p$-integer by Theorem \ref{combthm'}. In view of
the inequality $\ord_p(rR+p^{\al}R')\ge\al$, we have
$$\ord_p(R)\gs \al-\ord_p(r)=\tau_p(\{r\}_{p^{\al}},\{n-r\}_{p^{\al}}),$$
where the last equality follows from the definition of $\tau_p$
and the condition $n\eq0\not\eq r\ (\mo\ p^{\al})$.

{\it Case 2}. $\ord_p(r)\ls\beta<\al$. Since $\lfloor
n/p^{\al}\rfloor=\lfloor(n-1)/p^{\al}\rfloor$, the definition of
$R'$ implies that
$$\f{p^\al R'}n=\f1{\lfloor (n-1)/p^{\al}\rfloor!}
\sum_{k\eq r-1\,(\mo\ p^{\al})}\bi
{n-1}k(-1)^k\l(\f{k-(r-1)}{p^{\al}}\r)^{l}.$$ Applying the
induction hypothesis, we find that
$$\ord_p(p^{\al}R')-\beta\ge\tau_p(\{r-1\}_{p^{\al}},\{n-1-(r-1)\}_{p^{\al}})
=\tau_p(\{r-1\}_{p^{\al}},\{n-r\}_{p^{\al}}).$$ Since
$\{r\}_{p^{\al}}+\{n-r\}_{p^{\al}}\eq n\not\eq0\ (\mo\ p^{\al})$
and
\begin{eqnarray*}\bi{\{r\}_{p^{\al}}+\{n-r\}_{p^{\al}}}{\{r\}_{p^{\al}}}
&=&\f{\{r\}_{p^{\al}}+\{n-r\}_{p^{\al}}}{\{r\}_{p^{\al}}}
\bi{\{r\}_{p^{\al}}+\{n-r\}_{p^{\al}}-1}{\{r\}_{p^{\al}}-1}
\\&=&\f{\{r\}_{p^{\al}}+\{n-r\}_{p^{\al}}}{\{r\}_{p^{\al}}}
\bi{\{r-1\}_{p^{\al}}+\{n-r\}_{p^{\al}}}{\{r-1\}_{p^{\al}}},
\end{eqnarray*}
we have
$$\tau_p(\{r\}_{p^{\al}},\{n-r\}_{p^{\al}})
= \tau_p(\{r-1\}_{p^{\al}},\{n-r\}_{p^{\al}})+\beta-\ord_p(r).$$
Thus
$$\ord_p(p^{\al}R')\ge\ord_p(r)+\tau_p(\{r\}_{p^{\al}},\{n-r\}_{p^{\al}}).$$
Clearly
$\tau_p(\{r\}_{p^{\al}},\{n-r\}_{p^{\al}})\ls\al-\ord_p(r)$ by the
definition of $\tau_p$, so we also have
$$\ord_p(rR+p^{\al}R')\ge\ord_p(r)+\tau_p(\{r\}_{p^{\al}},\{n-r\}_{p^{\al}}).$$
Therefore
$$\ord_p(R)=\ord_p(rR)-\ord_p(r)\ge\tau_p(\{r\}_{p^{\al}},\{n-r\}_{p^{\al}}).$$

{\it Case 3}. $\beta<\min\{\al,\ord_p(r)\}$. In this case,
$\ord_p(\bar r)=\beta<\al$ where $\bar r=n-r$. Also,
\begin{eqnarray*}\sum_{k\eq \bar r\,(\mo\ p^{\al})}\bi nk(-1)^k\l(\f{k-\bar r}{p^{\al}}\r)^l
&=&\sum_{n-k\eq r\,(\mo\ p^{\al})}\bi
nk(-1)^k\l(\f{r-(n-k)}{p^{\al}}\r)^l
\\&=&(-1)^{l+n}\sum_{k\eq r\,(\mo\ p^{\al})}\bi nk(-1)^k\l(\f{k-r}{p^{\al}}\r)^l.
\end{eqnarray*}
Thus, as in the second case, we have
\begin{eqnarray*}\ord_p(R)&=&\ord_p\(\f1{\lfloor n/p^{\al}\rfloor!}
\sum_{k\eq\bar r\,(\mo\ p^{\al})} \bi nk(-1)^k\l(\f{k-\bar
r}{p^{\al}}\r)^l\)
\\&\ge&\tau_p(\{\bar r\}_{p^{\al}},\{n-\bar r\}_{p^{\al}})
=\tau_p(\{r\}_{p^{\al}},\{n-r\}_{p^{\al}}).
\end{eqnarray*}

The induction proof of Theorem \ref{strthm} is now complete.
\end{pf}

The following conjecture is based on extensive {\tt Maple}
calculations.
\begin{conj}\label{strconj} Let $p$ be any prime.
And let $\al,l\in\N$,  $n,r\in\Z$, with $n\ge2p^\a-1$. Then
equality in Theorem \ref{strthm} is attained if $l\ge\lfloor
n/p^\a\rfloor$ and
$$l\equiv\l\lfloor\f r{p^{\al}}\r\rfloor
+\l\lfloor\f{n-r}{p^{\al}}\r\rfloor\ \l(\mo\
(p-1)p^{\lfloor\log_p(n/p^\a)\rfloor}\r).$$
\end{conj}

\begin{rmk}\label{rmks} {\rm (1)
 The conjecture, if proved, would show that
 Theorem \ref{strthm} would be optimal in the sense that it is sharp
for infinitely many values of $l$.

(2) Note that the conjecture only deals with equality when
$l\ge\lfloor n/p^\a\rfloor$. For smaller values of $l$, our
inequality is still true, but not so strong. In \cite{SD}, we
obtain a stronger inequality when $l<\lfloor n/p^\a\rfloor$.}
\end{rmk}
\smallskip

We close with a table showing the amount by which the bound in
Theorem \ref{strthm} improves on that of Theorem \ref{combthm'}.
That is, we tabulate $\tau_p(\{r\}_{p^{\al}},\{n-r\}_{p^{\al}})$
when $p=3$ and $\a=2$.
\def\z{$0$}
\def\o{$1$}
\def\u{$2$}
\begin{table}[h]
\begin{center}
\caption{Values of $\tau_3(\{r\}_{9},\{n-r\}_{9})$}\label{comtbl}
\begin{tabular}{cc|ccccccccc|}
&\multicolumn{10}{c}{$\{r\}_9$}\\
&&\z&\o&\u&$3$&$4$&$5$&$6$&$7$&$8$\\
\cline{2-11}
&\z&\z&\u&\u&\o&\u&\u&\o&\u&\u\\
&\o&\z&\z&\u&\o&\o&\u&\o&\o&\u\\
&\u&\z&\z&\z&\o&\o&\o&\o&\o&\o\\
$\{n\}_9$&$3$&\z&\o&\o&\z&\u&\u&\o&\u&\u\\
&$4$&\z&\z&\o&\z&\z&\u&\o&\o&\u\\
&$5$&\z&\z&\z&\z&\z&\z&\o&\o&\o\\
&$6$&\z&\o&\o&\z&\o&\o&\z&\u&\u\\
&$7$&\z&\z&\o&\z&\z&\o&\z&\z&\u\\
&$8$&\z&\z&\z&\z&\z&\z&\z&\z&\z\\
\cline{2-11}
\end{tabular}
\end{center}
\end{table}

\def\line{\rule{.6in}{.6pt}}


\begin{thebibliography}{99}
\bibitem{BDSU} M. Bendersky and D. M. Davis, {\em 2-primary
$v_1$-periodic homotopy groups of ${\rm SU}(n)$}, Amer. J. Math.
{\bf 114} (1991) 529--544.
\bibitem{Bott} R. Bott, {\em The stable homotopy of the classical groups},
Annals of Math. {\bf 70} (1959) 313-337.
\bibitem{CMN1} F. R. Cohen, J. C. Moore, and
J. A. Neisendorfer, {\em The double suspension and exponents of
the homotopy groups of spheres}, Annals of Math. {\bf 110} (1979)
549--565.
\bibitem{CMN2} \line, {\em Exponents in homotopy theory}, Annals
of Math. Studies {\bf 113} (1987) 3--34.
\bibitem{large} D. M. Davis, {\em Elements of large order in
$\pi_*({\rm SU}(n))$}, Topology {\bf 37} (1998) 293--327.
\bibitem{DSU} \line, {\em $v_1$-periodic homotopy groups of ${\rm SU}(n)$ at
odd primes}, Proc. London Math. Soc. {\bf 43} (1991) 529--541.
\bibitem{sur} \line, {\em Computing $v_1$-periodic homotopy groups of spheres
and certain Lie groups}, Handbook of Algebraic Topology, Elsevier,
1995, pp. 993--1049.
\bibitem{DM} D. M. Davis and M. Mahowald, {\em Some remarks on
$v_1$-periodic homotopy groups}, London Math. Soc. Lect. Notes
{\bf 176} (1992) 55--72.
\bibitem{DP} D. M. Davis and K. Potocka, {\em 2-primary $v_1$-periodic
homotopy groups of ${\rm SU}(n)$ revisited}, to appear in Forum
Math., on-line version: {\tt
http://www.lehigh.edu/$\sim$dmd1/sun2long.pdf}.
\bibitem{Gray} B. Gray, {\em On the sphere of origin of infinite families in the
homotopy groups of spheres}, Topology {\bf 8} (1969) 219-232.
\bibitem{Hu} D. Husemoller, {\em Fibre Bundles}, 2nd ed., Springer, 1975.
\bibitem{IMJ} I. M. James, {\em On the suspension sequence}, Annals of Math. {\bf 65} (1957) 74-107.
\bibitem{LW} J.H. van Lint and R. M. Wilson, {\em A Course in Combinatorics}, 2nd ed.,
Cambridge Univ. Press, Cambridge, 2001.
\bibitem{Neis} J. A. Neisendorfer, {\em A survey of
Anick-Gray-Theriault constructions and applications to exponent
theory of spheres and Moore spaces}, Contemp. Math. Amer. Math.
Soc. {\bf 265} (2000) 159--174.
\bibitem{Selick} P. Selick, {\em Odd-primary torsion in $\pi_k(S^3)$}, Topology {\bf 17} (1978) 407-412.
\bibitem{S05} Z. W. Sun, {\em Polynomial extension of Fleck's congruence},
Acta Arith {\bf 122} (2006) 91-100.
\bibitem{SD} Z. W. Sun and D. M. Davis, {\em Combinatorial congruences modulo prime powers},
to appear in Trans. Amer. Math. Soc., on-line version: {\tt
http://arxiv.org/abs/math.NT/0508087}.
\bibitem{Th1} S. D. Theriault, {\em  2-primary exponent bounds for Lie groups of low rank},
Canad. Math. Bull. {\bf 47} (2004) 119--132.
\bibitem{Th2} \line, {\em The 5-primary homotopy exponent of the exceptional Lie group $E_8$},
J. Math. Kyoto Univ {\bf 44} (2004) 569--593.
\bibitem{To} H. Toda, {\em A topological proof of theorems of
Bott and Borel-Hirzebruch for homotopy groups of unitary groups},
Mem. Coll. Sci. Univ. Kyoto {\bf 32} (1959) 103-119.
\bibitem{W77} C. S. Weisman, {\em Some congruences for binomial coefficients},
Michigan Math. J. {\bf 24} (1977) 141--151.
\end{thebibliography}
\end{document}